\documentclass{article} 
\usepackage{iclr2024_workshop,times}

\usepackage{graphicx}
\usepackage{subcaption}
\usepackage{amsmath,amsfonts, amsthm, amssymb, bm}
\usepackage{authblk}
\usepackage{amsmath,amsfonts, amsthm, amssymb, bm}

\def\eqref#1{equation~\ref{#1}}

\def\1{\bm{1}}

\def\rvw{{\mathbf{w}}}
\def\rvx{{\mathbf{x}}}

\def\rmI{{\mathbf{I}}}

\def\rmX{{\mathbf{X}}}

\def\va{{\bm{a}}}
\def\vb{{\bm{b}}}

\def\vf{{\bm{f}}}

\def\vx{{\bm{x}}}

\def\mD{{\bm{D}}}

\DeclareMathAlphabet{\mathsfit}{\encodingdefault}{\sfdefault}{m}{sl}
\SetMathAlphabet{\mathsfit}{bold}{\encodingdefault}{\sfdefault}{bx}{n}

\def\gH{{\mathcal{H}}}

\newcommand{\E}{\mathbb{E}}

\newcommand{\R}{\mathbb{R}}

\newcommand{\norm}[1]{||#1||}
\newcommand{\loss}{\mathcal{E}}
\newcommand{\balpha}{\bm\alpha}
\DeclareMathOperator{\Dim}{dim}
\usepackage{hyperref}
\usepackage{url}

\title{Learning Stochastic Dynamics from Data}

\author[1]{Ziheng Guo}
\author[1]{Igor Cialenco}
\author[1]{Ming Zhong}
\affil[1]{Department of Applied Mathematics, Illinois Institute of Technology}
\begin{document}
\maketitle
\begin{abstract}
We present a noise guided trajectory based system identification method for inferring the dynamical structure from observation generated by stochastic differential equations.  Our method can handle various kinds of noise, including the case when the the components of the noise is correlated.  Our method can also learn both the noise level and drift term together from trajectory.  We present various numerical tests for showcasing the superior performance of our learning algorithm.
\end{abstract}
\section{Introduction}
Stochastic Differential Equation (SDE) is a fundamental modeling tool in various science and engineering fields~\cite{Evans, Sarkka_Solin_2019}. Compared with traditional deterministic models which often fall short in capturing the stochasticity in the nature, the significance of SDE lies in the ability to model complex systems influenced by random perturbations. Hence, SDE provides insights into the behavior of such systems under uncertainty.  By incorporating a random component, typically through a Brownian motion, SDE provides a more realistic and flexible framework for simulating and predicting the behavior of these complex and dynamic systems.

We consider the models following the SDE of the following form, $d\rvx_t = \vf(\rvx_t)dt + d\rvw_t$, where the state vector $\rvx_t \in \R^d$, the drift term $\vf: \R^d \rightarrow \R^d$, and the stochastic noise $\rvw_t$ is a Brownian motion with covariance matrix $\mD(\vx)$ which also depends on the state. In the field of mathematical finance, there are several SDE models that are widely used, for example, Black-Scholes Model \cite{0b9b8115-a8b8-3422-8e1c-a62077de6621} in option pricing, Vasicek Model \cite{VASICEK1977177} for analyzing interest rate dynamic and Heston Model \cite{10.1093/rfs/6.2.327} for volatility studies. Apart from the field of finance, SDE has its application in physics \cite{book, EGNS2008} and biology \cite{SZEKELY201414, DP2011} where SDE is used to study evolution of particles subjected to random forces and modeling physical and biological systems.

The calibration of these models is crucial for their effective application in areas like derivative pricing and asset allocation in finance, as well as in the analysis of particle systems in quantum mechanics and the study of biological system behaviors in biology. This requires the use of diverse statistical and mathematical techniques to ensure the models' outputs align with empirical data, thereby enhancing their predictive and explanatory power. Since these SDE models have explicit function form of both drift and diffusion terms, the calibration or estimation of parameters is mostly done by minimizing the least square error between observation and model prediction \cite{MrázekPospíšil+2017+679+704, 935092}. SDE has also been studied for more general cases where the drift term does not have an explicit form where it satisfies $\vf = \vf(\rvx_t, \theta)$ with $\theta$ being a vector a unknown parameters.  These models can be estimated with the maximum likelihood estimator approach presented in \cite{285f6c95-4239-301e-8a39-738a80e7080c}.  The other approach is to obtain the likelihood function by Radon Nikodym derivative over the whole observation $\{\rvx_t\}_{t=0}^{T}$.  Inspired by theorem $7.4$ in \cite{Sarkka_Solin_2019} and the Girsanov theorem, we derive a likelihood function that incorporates state-dependent correlated noise.  By utilizing such likelihood function, we are able to capture the essential structure of $\vf$ from data with complex noise structure. 

The remainder of the paper is structured as follows, section \ref{sec:Learning Framework} outlines the framework we use to learn the drift term and the noise level (if not state dependent), we demonstrates the effectiveness of our learning by testing it on various cases summarized in section \ref{sec:Examples} and additional examples presented in Appendix \ref{sec:more_exampels}, we conclude our paper in section \ref{sec:conclude} with a few pointers for ongoing and future developments.
\section{Learning Framework}\label{sec:Learning Framework}
We consider the following SDE
\begin{equation}\label{eq:original_model}
d\rvx_t = \vf(\rvx_t)dt + d\rvw_t, \quad \rvx_t, \rvw_t \in \R^d, 
\end{equation}
where $\vf: \R^d \rightarrow \R^d$ is a drift term, and $\rvw$ represents the Brownian noise with a symmetric positive definite covariance matrix $\mD = \mD(\vx): \R^d \rightarrow \R^{d \times d}$. 

We consider the scenario when we are given continuous observation data in the form of $\{\rvx_t, d\rvx_t\}_{t \in [0, T]}$ for $\rvx_0 \sim \mu_0$.  We will find the minimizer to the following loss function
\begin{equation}\label{eq:original_loss}
\loss_{\gH}(\tilde\vf) = \E_{\rvx_0 \in \mu_0}\Big[\frac{1}{2T}\int_{t = 0}^T<\tilde\vf(\rvx_t), \mD^{-1}(\rvx_t)\tilde\vf(\rvx_t)> \, dt - 2\int_{t = 0}^T<\tilde\vf(\rvx_t), \mD^{-1}(\rvx_t)d\rvx_t>\Big],
\end{equation}
for $\tilde\vf \in \gH$; the function space $\gH$ is designed to be convex and compact and it is also data-driven.  This loss functional is derived from the Girsanov theorem as well as inspiration from Theorem $7.4$ from \cite{Sarkka_Solin_2019}.

In the case of uncorrelated noise, i.e. $\mD(\vx) = \sigma^2(\vx)\rmI$, where $\rmI$ is the $d \times d$ identity matrix and $\sigma: \R^d \rightarrow \R^+$ is a scalar function depends on the state, representing the noise level, \ref{eq:original_loss} can be simplified to
\begin{equation}\label{eq:simple_loss}
\loss_{\gH}^{\text{Sim}}(\tilde\vf) = \E_{\rvx_0 \in \mu_0}\Big[\frac{1}{2T}\int_{t = 0}^T\frac{\norm{\tilde\vf(\rvx_t)}^2}{\sigma^2(\rvx_t)}\, dt - 2\int_{t = 0}^T\frac{<\tilde\vf(\rvx_t), d\rvx_t>}{\sigma^2(\rvx_t)}\Big].
\end{equation}
Moreover, we provide three different performance measures of our estimated drift.  First, if we have access to original drift function $\vf$, then we will use the following error to compute the difference between $\hat\vf$ (our estimator) to $\vf$ with the following norm
\begin{equation}\label{eq:rho_norm}
    \norm{\vf - \hat\vf}_{L^2(\rho)}^2 = \frac{1}{T}\int_{\rvx \in \Omega}\norm{\vf(\rvx) - \hat\vf(\rvx)}_{\ell^2(\R^d)}^2 \, d\rho(\rvx).
\end{equation}
Here the weighted measure $\rho$ is defined on $\Omega$, where it defines the region of $\rvx$ explored due to the dynamics defined by \eqref{eq:original_model}; therefore $\rho$ is given as follows
\begin{equation}\label{eq:rho_def}
\rho(\rvx) = \E_{\rvx_0 \sim \mu_0}\Big[\frac{1}{T}\int_{t = 0}^T\delta_{\rvx_t}(\rvx)\Big], \quad \text{where $\rvx_t$ evolves from $\rvx_0$ by \eqref{eq:original_model}.}
\end{equation}
The norm given by \eqref{eq:rho_norm} is useful only from the theoretical perspective, i.e. showing convergence.  In real life situation, $\vf$ is most likely non-accessible.  Hence we will look at a performance measure that compare the difference between $\rmX(\vf, \rvx_0, T) = \{\rvx_t\}_{t \in [0, T]}$ (the observed trajectory evolves from $\rvx_0 \sim \mu_0$ with the unknown $\vf$) and $\hat\rmX(\hat\vf, \rvx_0, T) = \{\hat\rvx_t\}_{t \in [0, T]}$ (the estimated trajectory evolves from the same $\rvx_0$ with the learned $\hat\vf$ and the same random noise as used by the original dynamics).  Then, the difference between the two trajectories is measured as follows
\begin{equation}\label{eq:traj_norm}
\norm{\rmX - \hat\rmX} = \E_{\rvx_0 \sim \mu_0}\Big[\frac{1}{T}\int_{t = 0}^T\norm{\rvx_t - \hat\rvx_t}_{\ell^2(\R^d)}^2 \, dt\Big].
\end{equation}
However, comparing two sets of trajectories (even with the same initial condition) on the same random noise is not realistic.  We compare the distribution of the trajectories over different initial conditions and all possible noise at some chosen time snapshots using the Wasserstein distance.
\section{Comparison to Others}
When the noise level becomes a constant, i.e. $\sigma(\rvx) = \sigma > 0$, we end up a much simpler loss
\[
\loss_{\gH}^{\text{Simpler}}(\tilde\vf) = \E_{\rvx_0 \in \mu_0}\Big[\frac{1}{2T\sigma^2}\int_{t = 0}^T\norm{\tilde\vf(\rvx_t)}^2\, dt - 2\int_{t = 0}^T<\tilde\vf(\rvx_t), d\rvx_t>\Big],
\]
which has been investigated in \cite{lu2022}.  System identification of the drift term has been studied in many different scenarios, e.g. identification by enforcing sparsity such as SINDy~\cite{sindy}, neural network based methods such as NeuralODE~\cite{neuralode2018} and PINN~\cite{pinn2019}, regression based~\cite{CS02}, and high-dimensional reduction framework~\cite{lu2019nonparametric}.  The uniqueness of our method is that we incorporate the covariance matrix into the learning and hence improving the estimation especially when the noise is correlated.
\section{Examples}\label{sec:Examples}
In this section, we demonstrate the application of our trajectory-based method for estimating drift functions, showcasing a variety of examples. We explore drift functions ranging from polynomials in one and two dimensions to combinations of polynomials and trigonometric functions, as well as a deep learning approach in one dimension. The observations, serving as the input dataset for testing our method, are generated by the Euler-Maruyama scheme, utilizing the drift functions as we just mentioned. The basis space $\gH$ is constructed employing either B-spline or piecewise polynomial methods for maximum degree p-max equals $2$. For higher order dimensions where $d \geq 2$, each basis function is derived through a tensor grid product, utilizing one-dimensional basis defined by knots that segment the domain in each dimension. The parameters for the following examples are listed in table \ref{tab: Parameters}. Due to space constraints, additional examples are provided in appendix \ref{sec:more_exampels}.

The estimation results are evaluated using several different metrics. We record the noise terms, $d\rvw_t$, from the trajectory generation process and compare the trajectories produced by the estimated drift functions, $\hat{\vf}$, under identical noise conditions. We examine trajectory-wise errors using equation \ref{eq:traj_norm} with relative trajectory error and plot both $\vf$ and $\hat{\vf}$ to calculate the relative $L^2$ error using \ref{eq:rho_norm}, where $\rho$ is derived by \ref{eq:rho_def}. Additionally, we assess the distribution-wise discrepancies between observed and estimated results, computing the Wasserstein distance at various time steps. 

\begin{table}[h!]
    \centering
\caption{Parameters Setup for Examples}
\label{tab: Parameters}
    \begin{tabular}{|c|c|c|c|} \hline 
         \multicolumn{2}{|c|}{Simulation Scheme}&  \multicolumn{2}{|c|}{Euler-Maruyama}\\ 
 \hline\hline
         $T$&  1&  $D$ $(d=1)$&0.6\\ \hline 
         $dt$&  0.001&  $D$ $(d=2)$& $\big(\begin{smallmatrix}
  0.6 & 0\\
  0 & 0.8
\end{smallmatrix}\big)$  \\ \hline 
         $\rvx_0$ &  Uniform(0,10)&  $M$& 10000\\ \hline 
         p-max&  2&  Basis Type& B-Spline / PW-Polynomial\\ \hline
    \end{tabular}
\end{table}
\subsection{Example: sine/cosine drift}
We initiate our numerical study with a one-dimensional ($d=1$) drift function that incorporates both polynomial and trigonometric components, given by $\vf = 2 + 0.08\rvx - 0.05\sin(\rvx) + 0.02\cos^2(\rvx)$. 

\begin{figure}[h!]
    \centering
    \includegraphics[width=0.75\linewidth]{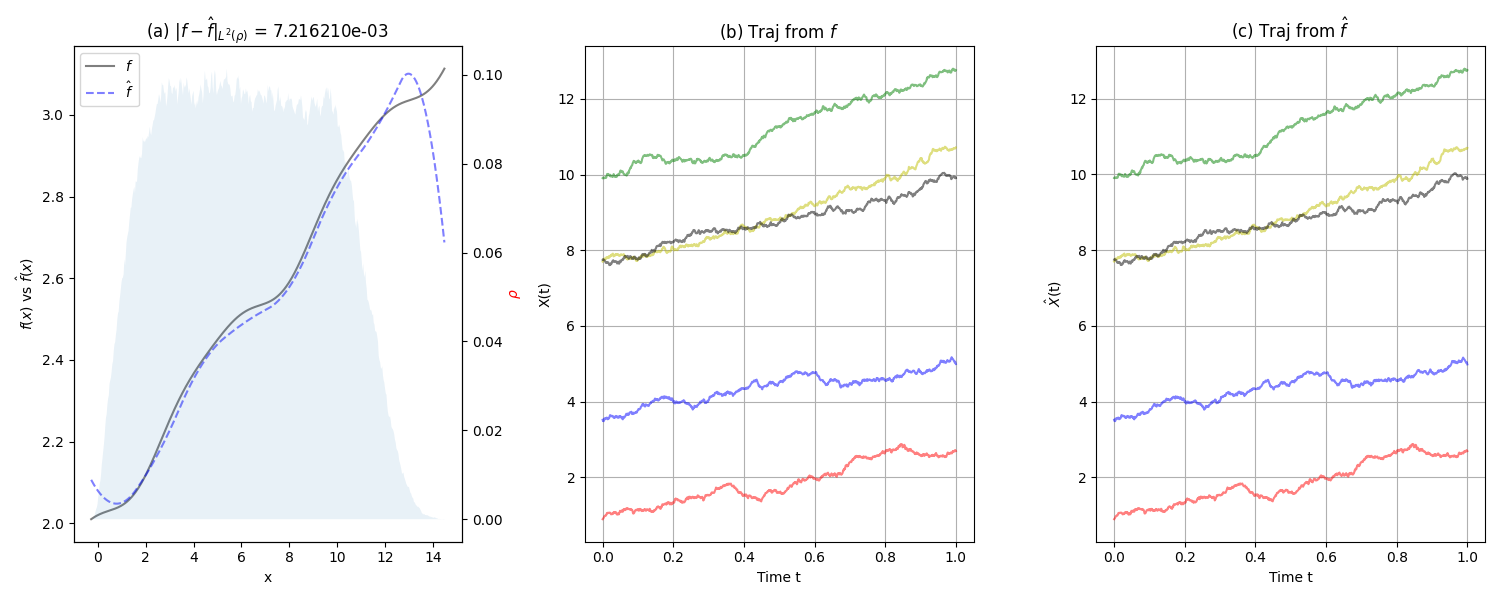}
    \caption{Left: Comparison of $\vf$ and $\hat{\vf}$. Middle: 5 trajectories generated by $\vf$. Right: 5 trajectories generated by $\hat{\vf}$ with same noise.}
    \label{fig:1D Example}
\end{figure}

Figure \ref{fig:1D Example} illustrates the comparison between the true drift function $\vf$ and the estimated drift function $\hat{\vf}$, alongside a comparison of trajectories. Notably, Figure \ref{fig:1D Example}(a) on the left includes a background region depicting the histogram of $\rvx_t$, which represents the distribution of observations over the domain of $\rvx$. This visualization reveals that in regions where $\rvx$ has a higher density of observations—indicated by higher histogram values—the estimation of $\hat{\vf}$ tends to be more accurate. Conversely, in less dense regions of the dataset (two ends of the domain), the estimation accuracy of $\hat{\vf}$ diminishes. 

Table \ref{Tab:1D Example} presents a detailed quantitative analysis of the estimation results, including the $L^2$ norm difference between $\vf$ and $\hat{\vf}$, as well as the trajectory error. Furthermore, the table compares the distributional distances between $\rvx_t$ and $\hat{\rvx}_t$ at selected time steps, with the Wasserstein distance results included.
\begin{table}[h!]
\caption{One-dimensional Drift Function Estimation Summary}
    \centering
    \begin{tabular}{|c|c||c|c|} \hline 
         Number of Basis&  8&  \multicolumn{2}{|c|}{Wasserstein Distance}\\ \hline 
         Maximum Degree&  2&  $t=0.25$& 0.0291\\ \hline 
         Relative $L^2(\rho)$ Error&  0.007935&  $t=0.50$& 0.0319\\ \hline 
         Relative Trajectory Error&  0.0020239 $\pm$ 0.002046&  $t=1.00$& 0.0403\\ \hline
    \end{tabular}
    \label{Tab:1D Example}
\end{table}
\subsection{Example: Deep Neural Networks}
\begin{figure}[h!]
    \centering
    \includegraphics[width=0.75\linewidth]{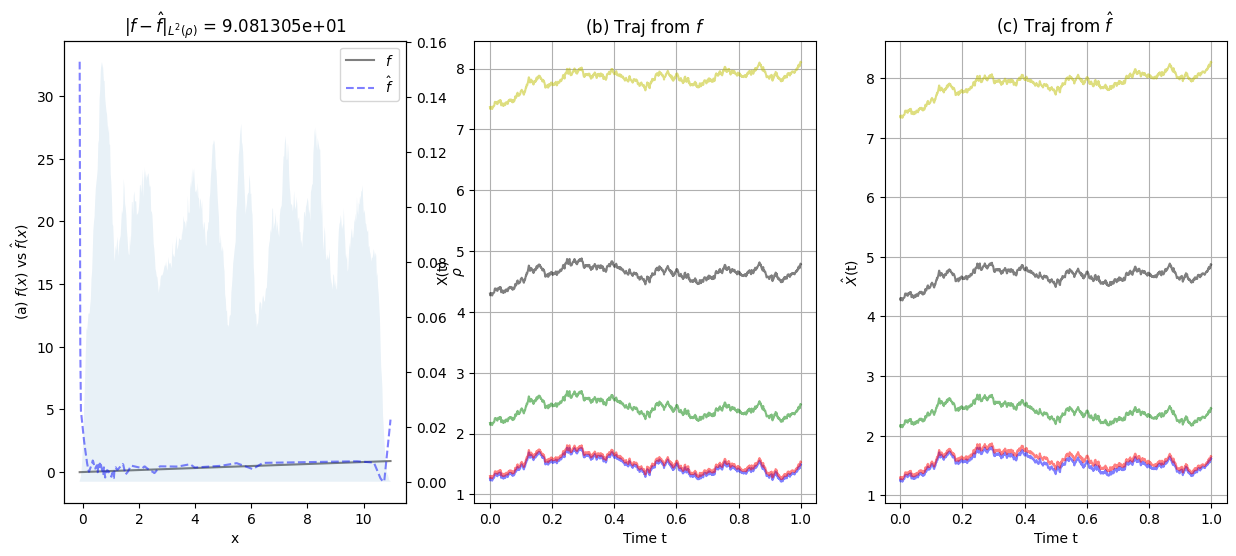}
    \caption{Left: Comparison of $\vf$ and $\hat{\vf}$. Middle: 5 trajectories generated by $\vf$. Right: 5 trajectories generated by $\hat{\vf}$ with same noise.}
    \label{fig:1D Example DL}
\end{figure}

We continue our numerical investigation with a one-dimensional ($d=1$) drift function which is given by $\vf = 0.08\rvx$.  Figure \ref{fig:1D Example DL} illustrates the comparison between the true drift function $\vf$ and the estimated drift function $\hat{\vf}$, alongside a comparison of trajectories.  Recall the setup of figures are similar to the ones presented in previous section.
The error for learning $\vf$ turns out to be bigger, especially towards the two end points of the interval.  However, the errors happen mostly during the two end points of the data interval, where the distribution of the data appears to be small, i.e. few data present in the learning.  We are able to recover most the trajectory.
\section{Conclusion}\label{sec:conclude}
We have presented a novel approach of learning the drift and diffusion from a noise-guide likelihood function.  Such learning can handle general noise structure as long as the covariance information of the noise is known.  Unknown noise structure when the covariance is a simple scalar is still learnable through quadratic variation.  However learning state-dependent covariance noise is ongoing.  Testing on high-dimensional drift term is ongoing.  Our loss only incorporates likelihood on the observation data.  It is possible to add various penalty terms to enhance the model performance, such as the sparsity assumption enforced in SINDy and other related deep learning methods such as NeuralODE and PINN. 
\subsubsection*{Acknowledgments}
ZG developed the algorithm, analyzed the data and implemented the software package.  ZG develops the theory with IC and MZ.  MZ designed the research.  All authors wrote the manuscript.  MZ is partially supported by NSF-$2225507$ and the startup fund provided by Illinois Tech.
\bibliography{iclr2024_workshop}
\bibliographystyle{iclr2024_workshop}

\clearpage
\appendix
\section{Implementation}\label{sec:code}
In this section, we will discuss in details how the algorithm is implemented for our learning framework presented in section \ref{sec:Learning Framework}.  In real applications, continuous data is rarely obtainable, we will have to deal with discrete observation data, i.e. $\{\rvx_l^m\}_{l, m = 1}^{L, M}$ with $\rvx_0^m$ being i.i.d sample from $\mu_0$, where $\rvx_l^m = \rvx^{(m)}(t_l)$ and $0 = t_1 < \cdots < t_L = T$.  We use a discretized version of \ref{eq:original_loss}, 
\begin{equation}\label{eq:original_loss_dis}
\begin{aligned}
\loss_{L, M, \gH}(\tilde\vf) &= \frac{1}{2TM}\sum_{l, m = 1}^{L - 1, M}\Big(<\tilde\vf(\rvx_l^m), \mD^{-1}(\rvx_l^m)\tilde\vf(\rvx_l^m)>(t_{l + 1} - t_l) \\
&\quad - 2<\tilde\vf(\rvx_l^m), \mD^{-1}(\rvx_l^m)(\rvx_{l + 1}^m - \rvx_l^m)>\Big),
\end{aligned}
\end{equation}
for $\tilde\vf \in \gH$.  Moreover, we also assume that $\gH$ is a finite-dimensional function space, i.e. $\Dim(\gH) = n < \infty$.  Then for any $\tilde\vf \in \gH$, $\tilde\vf(\vx) = \sum_{i = 1}^n\va_{i}\psi_{i}(\vx)$, where $\va_{i} \in \R^d$ is a constant vector coefficient and $\psi_{i}: \Omega \rightarrow \R$ is a basis of $\gH$ and the domain $\Omega$ is constructed by finding out the $\min/\max$ of the components of $\rvx_t \in \R^d$ for $t \in [0, T]$.  We consider two scenarios for constructing $\psi_i$, one is to use pre-determined basis such as piecewise polynomials, Clamped B-spline, Fourier basis, or a mixture of all of the aforementioned ones; the other is to use neural networks, where the basis functions are also trained from data.  Next, we can put the basis representation of $\tilde\vf$ back to \eqref{eq:original_loss_dis}, we obtain the following loss based on the coefficients
\begin{equation}\label{eq:original_loss_dis_coeff}
\begin{aligned}
\loss_{L, M, \gH}(\{\va_{\eta}\}_{i = 1}^n) &= \frac{1}{2TM}\sum_{l, m = 1}^{L - 1, M}\Big(\sum_{i = 1}^n\sum_{j = 1}^n<\va_i\psi_i(\rvx_l^m), \mD^{-1}(\rvx_l^m)\va_j\psi_j(\rvx_l^m)>(t_{l + 1} - t_l) \\
&\quad - 2\sum_{i = 1}^n<\va_i\psi_i(\rvx_l^m), D^{-1}(\rvx_l^m)(\rvx_{l + 1}^m - \rvx_l^m)>\Big),
\end{aligned}
\end{equation}
In the case of covariance matrix $\mD$ being a diagonal matrix, i.e.
\[
\mD(\vx) = \begin{bmatrix} \sigma_1^2(\vx) & 0 & \cdots & 0 \\ 0 & \sigma_2^2(\vx) & \cdots & 0 \\ \vdots & \vdots & \ddots & \vdots \\ 0 & 0 & \cdots & \sigma_d^2(\vx) \end{bmatrix} \in \R^{d \times d}, \quad \sigma_i > 0, \quad i = 1, \cdots, d.
\]
Then \eqref{eq:original_loss_dis_coeff} can be re-written as
\[
\begin{aligned}
\loss_{L, M, \gH}(\{\va_{\eta}\}_{i = 1}^n) &= \frac{1}{2TM}\sum_{l, m = 1}^{L - 1, M}\Big(\sum_{i = 1}^n\sum_{j = 1}^n\sum_{k = 1}^d\frac{(\va_i)_k(\va_j)_k}{\sigma_k^2(\rvx_l^m)}\psi_i(\rvx_l^m)\psi_j(\rvx_l^m)(t_{l + 1} - t_l) \\
&\quad - 2\sum_{i = 1}^n\sum_{k = 1}^d\frac{(\va_i)_k(\rvx_{l + 1}^m - \rvx_l^m)_k}{\sigma_k^2(\rvx_l^m)}\psi_i(\rvx_l^m)\Big),
\end{aligned}
\]
Here $(\vx)_k$ is the $k^{th}$ component of any vector $\vx \in \R^d$.  We define $\balpha_k = \begin{bmatrix}(\va_1)_k & \cdots & (\va_n)_k \end{bmatrix}^\top \in \R^n$, and $A_k \in \R^{n \times n}$ as 
\[
A_k(i, j) = \frac{1}{2TM}\sum_{l, m = 1}^{L - 1, M}\Big(\sum_{i = 1}^n\sum_{j = 1}^n\frac{(\va_i)_k(\va_j)_k}{\sigma_k^2(\rvx_l^m)}\psi_i(\rvx_l^m)\psi_j(\rvx_l^m)(t_{l + 1} - t_l),
\]
and $\vb_k \in \R^{n}$ as
\[
\vb_k(i) = \sum_{i = 1}^n\frac{(\va_i)_k(\rvx_{l + 1}^m - \rvx_l^m)_k}{\sigma_k^2(\rvx_l^m)}\psi_i(\rvx_l^m)\Big).
\]
Then \eqref{eq:original_loss_dis_coeff} can be re-written as
\[
\loss_{L, M, \gH}(\{\va_{\eta}\}_{i = 1}^n) = \sum_{k = 1}^d(\balpha_k^\top A_k\balpha_k - 2\balpha_k^\top\vb_k).
\]
Since each $\balpha_k^\top A_k\balpha_k - 2\balpha_k^\top\vb_k$ is decoupled from each other, we just need to solve simultaneously 
\[
A_k\hat\balpha_k - \vb_k = 0, \quad k = 1, \cdots, d.
\]
Then we can obtain $\hat\vf(\vx) = \sum_{i = 1}^n\hat\va_i\psi_k(\vx)$.  However when $\mD$ does not have a diagonal structure, we will have to resolve to gradient descent methods to minimize \eqref{eq:original_loss_dis_coeff} in order to find the coefficients $\{\va_i\}_{i = 1}^n$ for a total number of $nd$ parameters.  However, if a data-driven basis is desired, then we set $\gH$ to be the space neural networks with the same depth, same number of neurons and same activation functions on the hidden layers.  Furthermore, we find $\hat\vf$ from minimizing \eqref{eq:original_loss_dis} using any deep learning optimizer such as Stochastic Gradient Descent or Adam from well-known deep learning packages. 
\section{Additional Examples}\label{sec:more_exampels}
\begin{figure}[h!]
    \centering
    \includegraphics[width=0.75\linewidth]{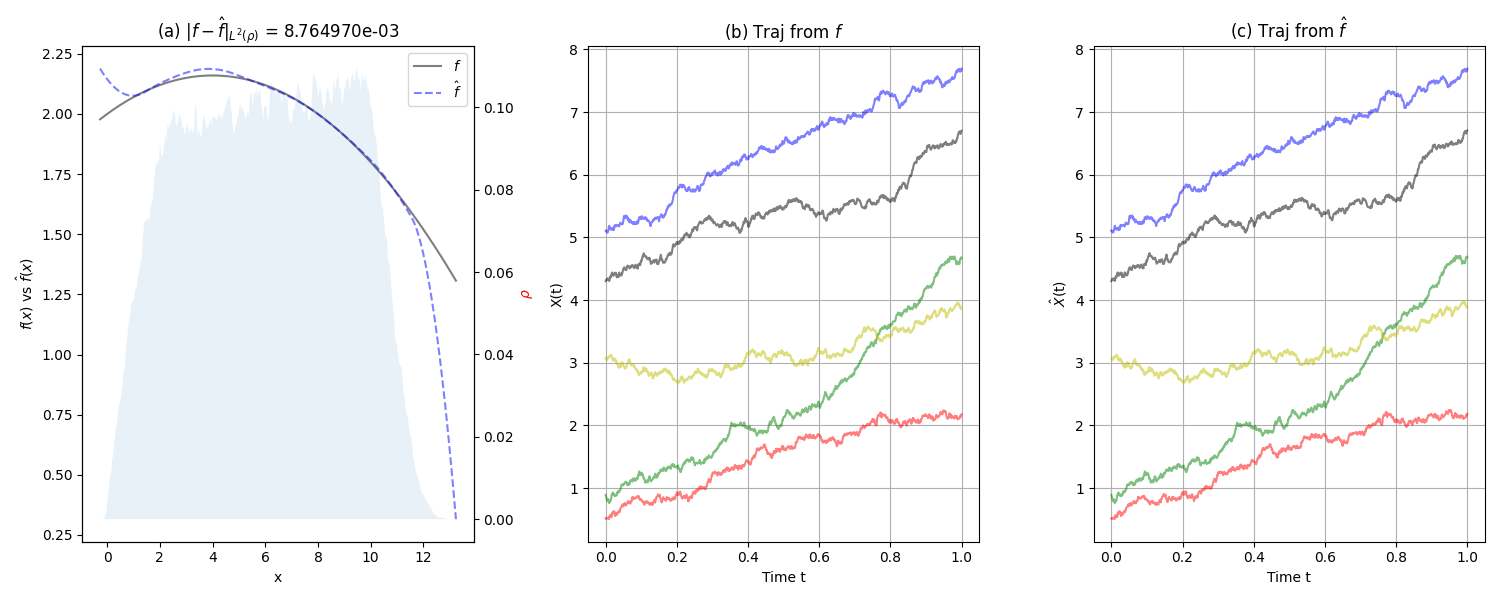}
    \caption{One-dimensional Polynomial Drift Comparison Summary}
    \label{fig:1D Poly Summary}
\end{figure}

For $d = 1$, we also worked on polynomial drift function $\vf = 2 + 0.08\rvx - 0.01\rvx^2$. The estimation results are depicted in \ref{fig:1D Poly Summary} and detailed in Table \ref{tab: 1D Poly Summary}. 

\begin{table}[h!]
\caption{One-dimensional Polynomial Drift Function Estimation Summary}
    \centering
    \begin{tabular}{|c|c||c|c|} \hline 
         Number of Basis&  10&  \multicolumn{2}{|c|}{Wasserstein Distance}\\ \hline 
         Maximum Degree&  2&  $t=0.25$& 0.0153\\ \hline 
         Relative $L^2(\rho)$ Error&  0.0087649&  $t=0.50$& 0.0154\\ \hline 
         Relative Trajectory Error&  0.00199719 $\pm$ 0.00682781&  $t=1.00$& 0.0278\\ \hline
    \end{tabular}
    \label{tab: 1D Poly Summary}
\end{table}

For $d = 2$, we examine two types of drift function $\vf$: polynomial and trigonometric. Denote
\begin{equation*}
    \vf = \begin{pmatrix}
  \vf_1 \\ 
  \vf_2
  \end{pmatrix}
    \quad\text{and}\quad
  \rvx =
  \begin{pmatrix}
  \rvx_1\\ 
  \rvx_2
  \end{pmatrix}
\end{equation*}
where $\vf_i: \R^2 \rightarrow \R$ and $\rvx_i \in \R$ for $i \in \{1, 2\}$.

\begin{figure}
    \centering
    \includegraphics[width=0.75\linewidth]{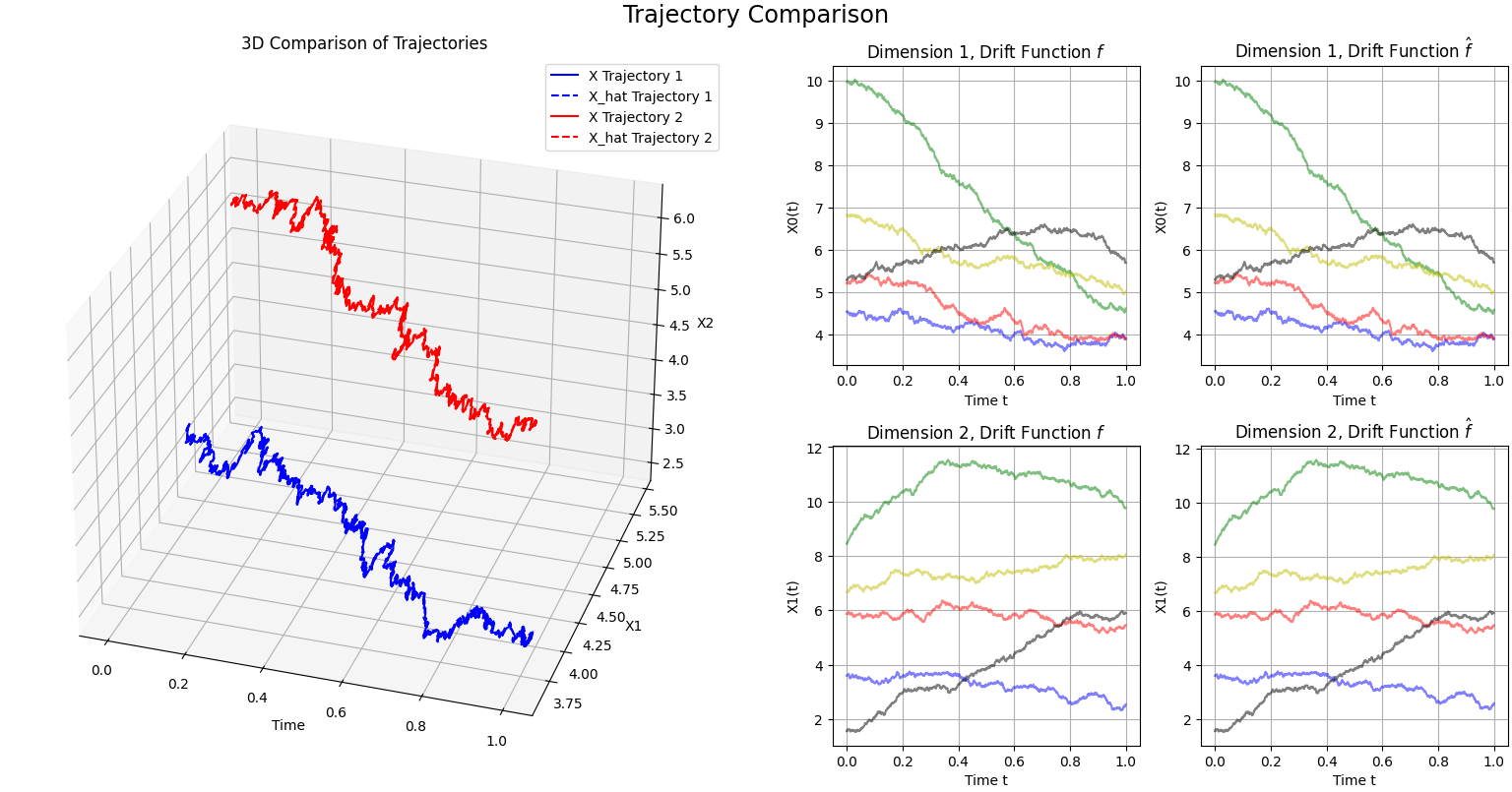}
    \caption{Two-dimensional Polynomial Trajectory Comparison}
    \label{fig:2D Traj}
\end{figure}

\begin{figure}
\centering
\begin{subfigure}[b]{0.75\textwidth}
   \includegraphics[width=1\linewidth]{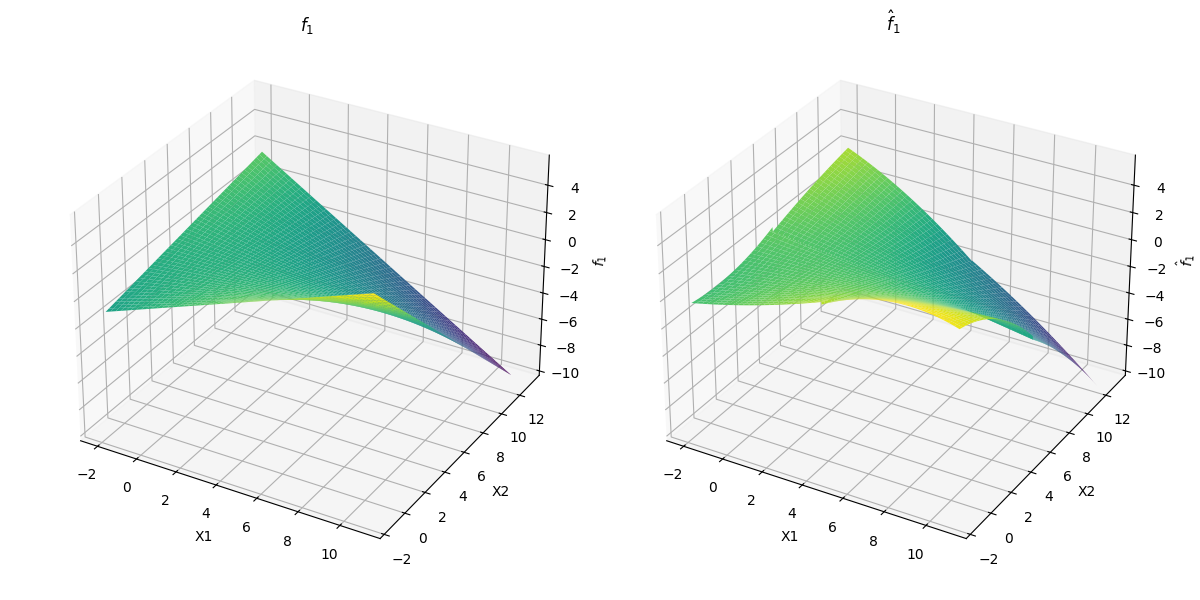}
   \caption{}
   \label{fig:2d poly Example 1} 
\end{subfigure}

\begin{subfigure}[b]{0.75\textwidth}
   \includegraphics[width=1\linewidth]{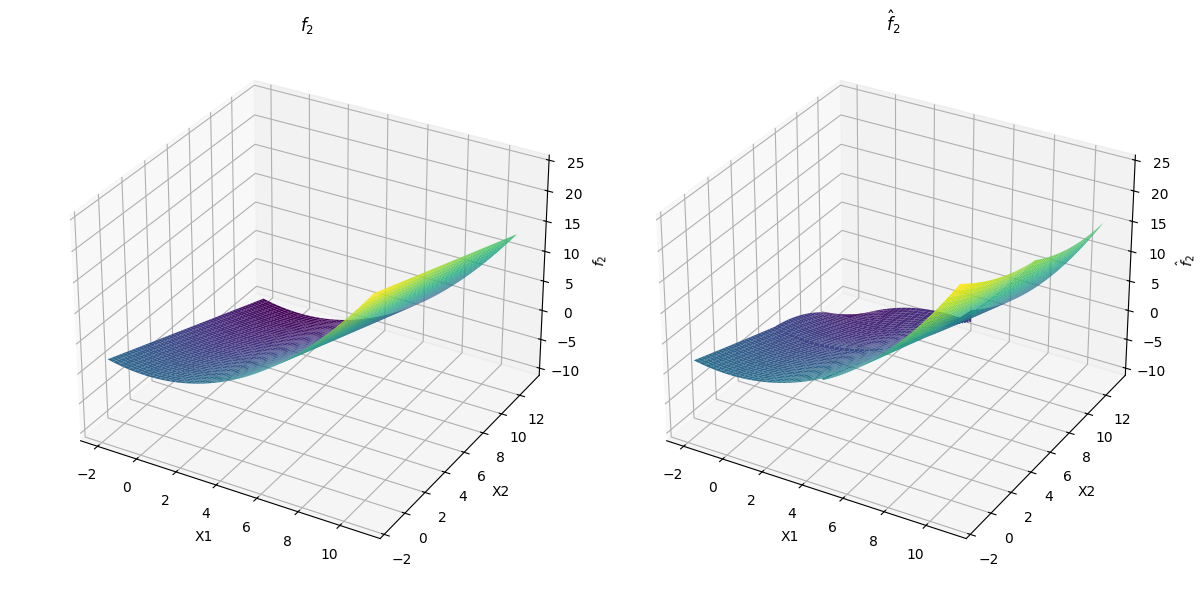}
   \caption{}
   \label{fig:2d poly Example 2}
\end{subfigure}

\caption[Comparison of $\vf$ and $\hat{\vf}$ in 2D]{Two-dimensional Polynomial Comparison of $\vf$ and $\hat{\vf}$. (a) Surface of Dimension 1 (b) Surface of Dimension 2}
\end{figure}

For polynomial drift function $\vf$, we set 
\begin{align*} 
    \vf_1 &= 0.4\rvx_1 - 0.1\rvx_1\rvx_2 \\
    \vf_2 &= -0.8\rvx_2 + 0.2\rvx_1^2.
\end{align*}
Figure \ref{fig:2D Traj}, Figure \ref{fig:2d poly Example 1}, \ref{fig:2d poly Example 2} and Table \ref{tab: 2D Poly} shows evaluation of the polynomial drift function estimation result. 

\begin{figure}
    \centering
    \includegraphics[width=0.75\linewidth]{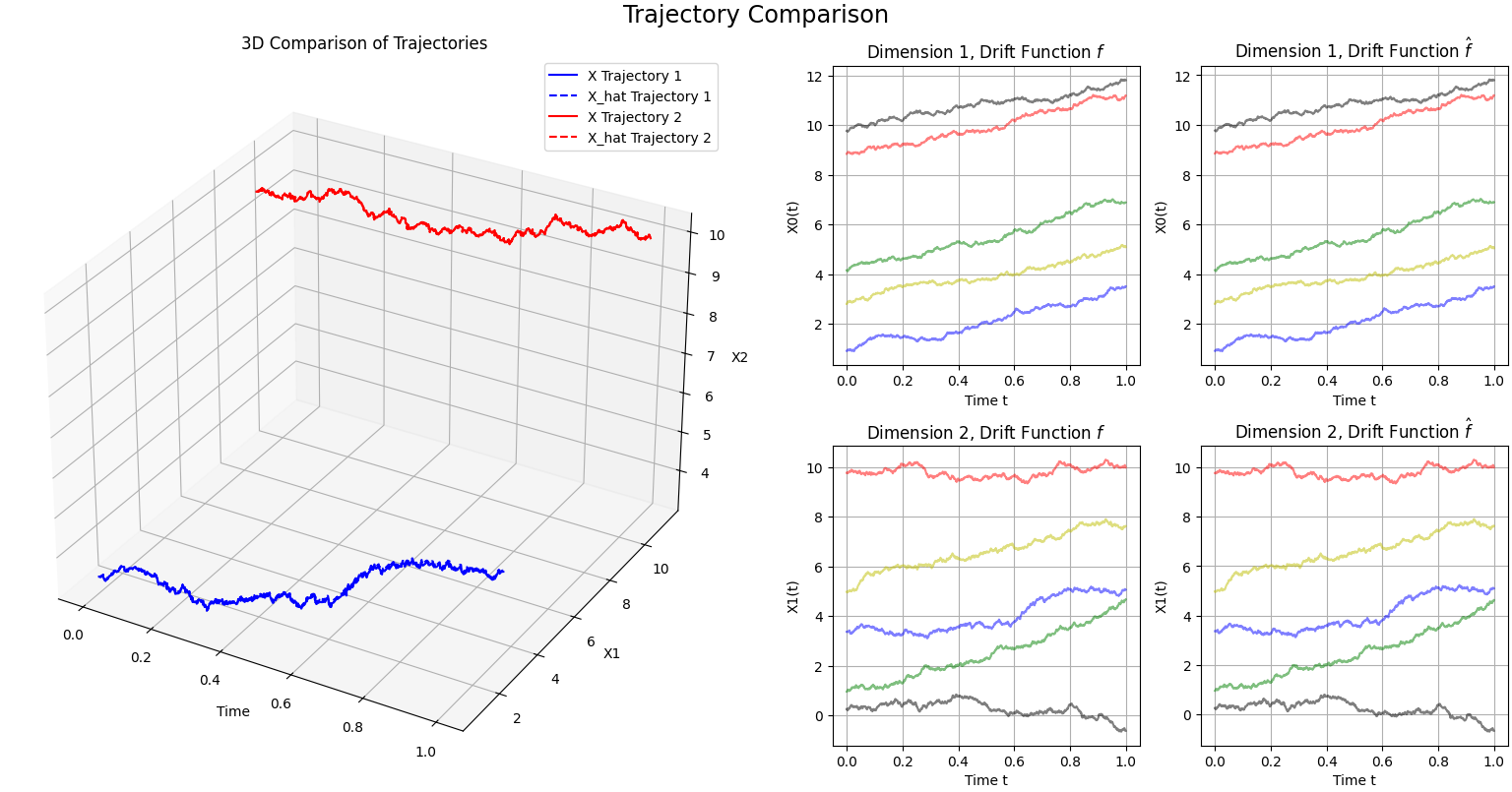}
    \caption{Two-dimensional Trigonometric Trajectory Comparison}
    \label{fig:Tri Traj}
\end{figure}

\begin{figure}
\centering
\begin{subfigure}[b]{0.75\textwidth}
   \includegraphics[width=1\linewidth]{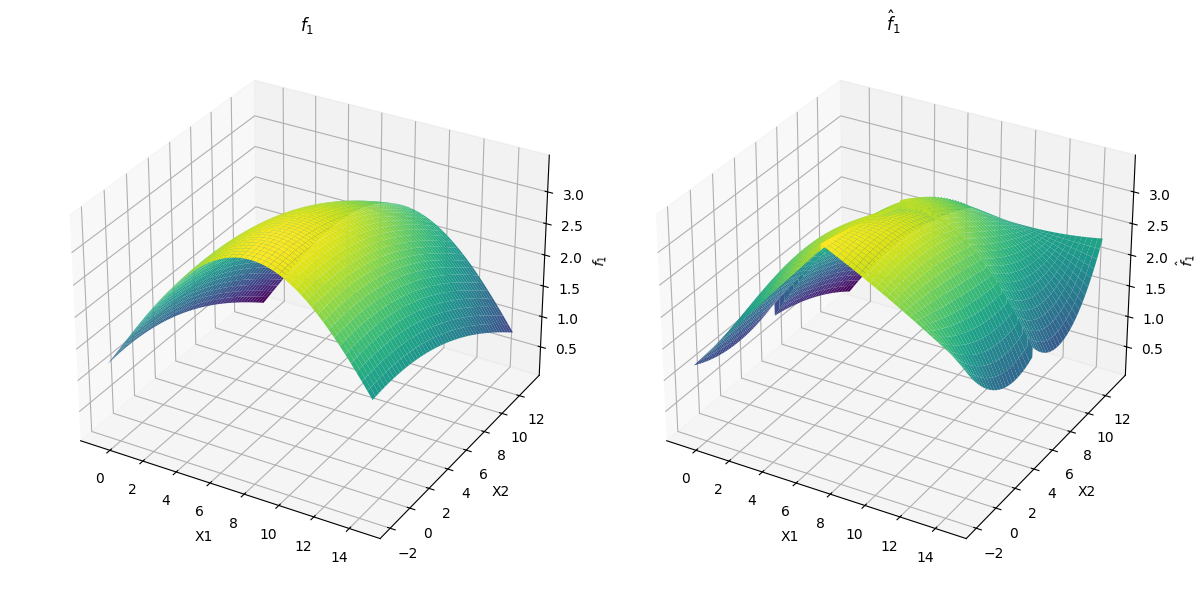}
   \caption{}
   \label{fig:2D Tri Example 1} 
\end{subfigure}

\begin{subfigure}[b]{0.75\textwidth}
   \includegraphics[width=1\linewidth]{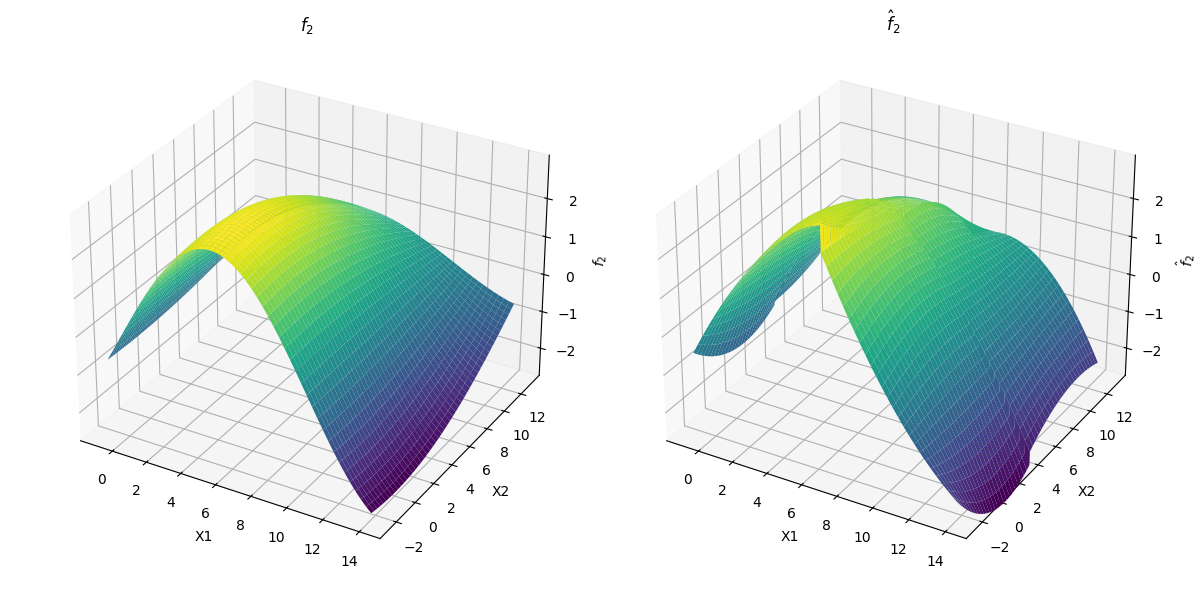}
   \caption{}
   \label{fig:2D Tri Example 2}
\end{subfigure}

\caption[Comparison of $\vf$ and $\hat{\vf}$ in 2D]{Two-dimensional Trigonometric Comparison of $\vf$ and $\hat{\vf}$. (a) Surface of Dimension 1 (b) Surface of Dimension 2}
\end{figure}

For trigonometric drift function $\vf$, we set 
\begin{align*} 
    \vf_1 &= 2\sin(0.2\rvx_1) + 1.5\cos(0.1\rvx_2) \\
    \vf_2 &=3\sin(0.3\rvx_1)\cos(0.1\rvx_2).
\end{align*}
Figure \ref{fig:Tri Traj}, Figure \ref{fig:2D Tri Example 1}, \ref{fig:2D Tri Example 2} and Table \ref{tab: 2D Tri} shows evaluation of the trigonometric drift function estimation result.

\begin{table}[h!]
\caption{Two-dimensional Polynomial Drift Function Estimation Summary}
    \centering
    \begin{tabular}{|c|c||c|c|} \hline 
         Number of Basis&  36&  \multicolumn{2}{|c|}{Wasserstein Distance}\\ \hline 
         Maximum Degree&  2&  $t=0.25$& 0.0891\\ \hline 
         Relative $L^2(\rho)$ Error&  0.02118531&  $t=0.50$& 0.0872\\ \hline 
         Relative Trajectory Error&  0.00306613 $\pm$ 0.00375144&  $t=1.00$& 0.0853\\ \hline
    \end{tabular}
    \label{tab: 2D Poly}
\end{table}

\begin{table}[h!]
\caption{Two-dimensional Trigonometric Drift Function Estimation Summary}
    \centering
    \begin{tabular}{|c|c||c|c|} \hline 
         Number of Basis&  36&  \multicolumn{2}{|c|}{Wasserstein Distance}\\ \hline 
         Maximum Degree&  2&  $t=0.25$& 0.1011\\ \hline 
         Relative $L^2(\rho)$ Error&  0.02734505&  $t=0.50$& 0.1119\\ \hline 
         Relative Trajectory Error&  0.0041613 $\pm$ 0.0079917&  $t=1.00$& 0.1293\\ \hline
    \end{tabular}
    \label{tab: 2D Tri}
\end{table}

\end{document}